\input amstex
\documentstyle{amsppt}
\document

\magnification 1100

\def\gen{{\frak{g}}}
\def\len{{\frak{l}}}

\def\ten{\frak{t}}

\def\uen{\frak{u}}

\def\sen{\frak{s}}
\def\pen{\frak{p}}

\def\qen{\frak{q}}
\def\Pen{\frak{P}}
\def\Cen{\frak{C}}
\def\Sen{\frak{S}}

\def\1b{{\bold 1}}
\def\ab{{\bold a}}
\def\bb{{\bold b}}

\def\eb{{\bold e}}
\def\fb{{\bold f}}

\def\Ab{{\bold A}}
\def\Bb{{\bold B}}

\def\Kb{{\bold K}}
\def\Hb{{\bold H}}

\def\Qb{{\bold Q}}

\def\Sb{{\bold S}}
\def\Tb{{\bold T}}
\def\Lb{{\bold L}}
\def\Ub{{\bold U}}

\def\Zb{{\bold Z}}

\def\b{{\roman b}}

\def\k{{\roman k}}

\def\A{{\roman A}}
\def\B{{\roman B}}

\def\K{{\roman K}}

\def\and{{\quad\text{and}\quad}}

\def\Im{\text{Im}\,}

\def\mod{\text{mod}\,}

\def\CC{{\Bbb C}}

\def\NN{{\Bbb N}}

\def\ZZ{{\Bbb Z}}

\def\Hc{{\Cal H}}

\def\eps{{\epsilon}}

\def\st{{\,|\,}}

\def\ds{\displaystyle}                
\def\ts{\textstyle}                
                
\def\qed{\hfill $\sqcap \hskip-6.5pt \sqcup$}        
\overfullrule=0pt                                    

\def\Umod{{\bold{\dot{U}}}}
\def\Bmod{{\bold{\dot{B}}}}
\def\Lmod{{\bold{\dot{L}}}}
\def\var{\varepsilon}

\newdimen\Squaresize\Squaresize=14pt
\newdimen\Thickness\Thickness=0.5pt
\def\Square#1{\hbox{\vrule width\Thickness
	      \vbox to \Squaresize{\hrule height \Thickness\vss
	      \hbox to \Squaresize{\hss#1\hss}
	      \vss\hrule height\Thickness}
	      \unskip\vrule width \Thickness}
	      \kern-\Thickness}
\def\Vsquare#1{\vbox{\Square{$#1$}}\kern-\Thickness}


\centerline{\bf GEOMETRIC CONSTRUCTION OF THE GLOBAL BASE}
\centerline{\bf OF THE QUANTUM MODIFIED ALGEBRA OF $\widehat{\gen\len}_N$}

\vskip1cm

\centerline{\bf O. SCHIFFMANN and E. VASSEROT
\footnote{
The second author is partially supported by EEC grant
no. ERB FMRX-CT97-0100.\hfill\break}}

\vskip3cm

\centerline{\smc 1. Introduction and notations.}

\vskip3mm

\noindent{\bf 1.1.}
A geometric construction of the modified quantum algebra of $\gen\len_n$
was given in \cite{BLM}. It was then observed
independentely by Lusztig and Ginzburg-Vasserot (see \cite{L1}, \cite{GV})
that this construction admits an affine analogue in terms of periodic flags
of lattices. However the compatibility of the canonical base of
the modified algebra and of the geometric base given by intersection cohomology
sheaves on the affine flag variety was never proved. The aim of the paper is 
to prove this compatibility.
As a consequence we prove a recent conjecture of Lusztig (see \cite{L1}). 
Of course, our proof would work also in the finite type case.

\vskip3mm

\noindent{\bf 1.2.} Let $v,z,$ be formal variables.
Set $\Zb=\CC[v]$, $\Ab=\CC[v,v^{-1}]$ and $\Kb=\CC(v)$.
Let $\k$ be the field with $q^2$ elements, where $q$ is a prime power.
Consider the local field $\K=\k((z))$ and the ring of integers 
$\A=\k[[z]]\subset\K$.
For any set $X$ with the action of a group $G$,
let $\CC_G[X]$ be the space of $G$-invariant complex functions on $X$
supported on a finite number of orbits.
Given $G$-sets $X_1, X_2, X_3$, we consider the convolution product
$$\CC_G[X_1\times X_2]\times\CC_G[X_2\times X_3]\to\CC_G[X_1\times X_3],\qquad
(\alpha,\beta)\mapsto(p_{13})_!\bigl((p_{12}^*\alpha)(p_{23}^*\beta)\bigr),$$
whenever it is well-defined.
Here $p_{ij}\,:\,X_1\times X_2\times X_3\to X_i\times X_j$ is the projection
along the factor not named.
If $X$ is an irreducible algebraic variety, let $\Hc^i(IC_X)$ be the
$i$-th cohomology sheaf of the intersection cohomology complex of $X$.
Then, for any stratum $Y\subseteq X$, let $\dim\Hc^i_Y(IC_X)$ be the dimension
of the stalk of $\Hc^i(IC_X)$ at a point in $Y$.
For any subset $Y\subset X$ let $\bar Y$ denote the Zariski closure 
of $Y$ in $X$. 

\vskip1cm

\centerline{\smc 2. Reminder on flag varieties.}

\vskip3mm

\noindent{\bf 2.1.}
Fix a positive integer $D$. Set $G_D=GL_D(\K)$.
A lattice in $\K^D$ is a free $\A$-submodule of rank $D$.
Let $Y_D$ be the set of {\bf $n$-steps periodic flags} in $\K^D$,
i.e. the set of sequences of lattices $L=(L_i)_{i\in\ZZ}$ such that
$L_i\subseteq L_{i+1}$ and $L_{i+n}=z^{-1}\,L_i.$
The group $G_D$ acts on $Y_D$ in the obvious way.
Let $\Sen_D$ be the set of all $\ZZ\times\ZZ$-matrices with
entries in $\NN$, say $\sen=(\sen_{ij})_{i,j\in\ZZ}$, such that 
$\sen_{i+n,j+n}=\sen_{ij}$ and
$${\ts\sum_{i\in\ZZ}\sum_{j=1}^n\sen_{ij}=D.}$$
The set $\Sen_D$ parametrizes the orbits of the diagonal action of
$G_D$ in $Y_D\times Y_D$ : to $\sen$ corresponds the set $Y_\sen$ of
the couple $(L,L')$ such that
$$\sen_{ij}=\dim\biggl({L_i\cap L'_j\over
(L_{i-1}\cap L'_j)+(L_i\cap L'_{j-1})}\biggr)$$
(see \cite{L1, Lemma 1.5} for instance).
For all $L\in Y_D$, $\sen\in\Sen_D$, the fiber over $L$ of the first 
projection $Y_\sen\to Y_D$ is either empty or the set of 
$\k$-points of an algebraic variety whose isomorphism class
is independent of $L$. Let $y_\sen$ be its dimension.
We have $y_\sen=\sum\sen_{ij}\sen_{kl},$
where the sum is taken over the set
$$\{(i,j,k,l)\st i\geq k,\,j<l,\,i\in[1,n]\}$$ 
(see \cite{L1, Lemma 4.3} for instance).
Let $\1b_\sen\in\CC_{G_D}[Y_D\times Y_D]$ be 
the characteristic function of $Y_\sen$.
The convolution product is well-defined and endows $\CC_{G_D}[Y_D\times Y_D]$
with the structure of an associative algebra. 

\vskip3mm

\noindent{\bf 2.2.}
Let $X_D$ be the set of {\bf complete periodic flags} in $\K^D$,
i.e. the set of sequences of lattices $L=(L_i)_{i\in\ZZ}$ such that
$L_i\subsetneq L_{i+1}$ and $L_{i+D}=z^{-1}\,L_i.$
The group $G_D$ acts on $X$ in the obvious way.
Let $\Pen_D$ be the set of functions $\pen\,:\,\ZZ\to\ZZ$ such that
$\pen(j+D)=\pen(j)+n$ for all $j$. 
The set $\Pen_D$ parametrizes the orbits of the diagonal action of
$G_D$ in $Y_D\times X_D$ : to $\pen\in\Pen_D$ corresponds the orbit $X_\pen$ of
the couple $(L_\pen,L_\emptyset)$ such that
$$L_{\pen,i}=\prod_{\pen(j)\leq i}\k\, e_j
\quad\text{and}\quad
L_{\emptyset,i}=\prod_{j\leq i}\k\, e_j,$$
where $(e_1,e_2,...,e_D)$ is a fixed $\K$-basis of
$\K^D$, and $e_{i+kD}=z^{-k}\, e_i$ for all $k\in\ZZ$.
If $L\in Y_D$, the fiber over $L$ of the 
projection $X_\pen\to Y_D$ is either empty or the set of 
$\k$-points of an algebraic variety whose isomorphism class is independent of 
$L$. Let $x_\pen$ be its dimension and let $\1b_\pen\in\CC_{G_D}[Y_D\times X_D]$
be the characteristic function of $X_\pen$.
We have 
$$x_\pen=\sharp\{(k,l)\st\pen(k)\in[1,n],\,k<l,\,\pen(k)\geq\pen(l)\}.$$ 
The convolution product is well-defined and endows the
space $\CC_{G_D}[Y_D\times X_D]$ with the structure of a left 
$\CC_{G_D}[Y_D\times Y_D]$-module
and a right $\CC_{G_D}[X_D\times X_D]$-module.

\vskip3mm

\noindent{\bf 2.3.}
Let $S^{^f}_D$ be the symmetric group and let
$\Hb_D$, $S_D$, be the affine Hecke algebra and the affine 
symmetric group of type $GL_D$. Recall that $\Hb_D$ is the 
unital associative $\Ab$-algebra generated by $T_i^{\pm 1},X_j^{\pm 1},$
$i\in[1,D-1],$ $j\in[1,D]$, modulo the following relations 
$$\matrix
T_i\,T_i^{-1}=1=T_i^{-1}\,T_i,\hfill\qquad&(T_i+1)(T_i-v^{-2})=0,\hfill\cr\cr
T_i\,T_{i+1}\,T_i=T_{i+1}\,T_i\,T_{i+1},\hfill\qquad&
|i-j|>1\Rightarrow T_i\,T_j=T_j\,T_i,\hfill\cr\cr
X_i\,X_i^{-1}=1=X_i^{-1}\,X_i,\hfill\qquad&X_i\,X_j=X_j\,X_i,\hfill\cr\cr
T_i\,X_i\,T_i=v^{-2}X_{i+1},\hfill\qquad& 
j\not= i,i+1\Rightarrow X_j\,T_i=T_i\,X_j.\hfill
\endmatrix$$
An element $\pen\in\Pen_D$ is identified with the $D$-uple
$\bigl(\pen(1),\pen(2),...,\pen(D)\bigr)\in\ZZ^D$. Let $\Cen_D\subset\Pen_D$
be the subset formed by the $D$-uples 
$\lambda=(\lambda_1,\lambda_2,...,\lambda_D)$ such that
$$1\leq\lambda_1\leq\lambda_2\leq...\leq\lambda_D\leq n.$$
The affine symmetric group $S_D$ acts on $\Pen_D$ on the right in the 
obvious way : elements of $S_D^{^f}$ act by permutations of the components on 
$\ZZ^D$ and
$$(\pen)\mu=(\pen(1)+n\mu_1,\pen(2)+n\mu_2,...,\pen(D)+n\mu_D),
\qquad\forall\mu\in\ZZ^D.$$
The set $\Cen_D$ is a fundamental domain for this action.
For any $\lambda\in\Cen_D$ let $S_\lambda$ be the Young subgroup
$$S_\lambda=\{w\in S_D\st (\lambda)w=\lambda\}\subseteq S^{^f}_D.$$
If $\lambda\in\Cen_D$, let $\Pen_\lambda$ be the $S_D$-orbit of $\lambda$
in $\Pen_D$. An element $\pen\in\Pen_\lambda$ is identified with a class
in $S_\lambda\setminus S_D$. Put $T_\pen=\sum_{w\in\pen}T_w$ for
all $\pen\in\Pen_D$, and set
$${\ts\Tb_D=\bigoplus_{\lambda\in\Cen_D}T_\lambda\Hb_D.}$$
The space $T_\lambda\Hb_D$ is linearly spanned
by the elements $T_\pen$ with $\pen\in\Pen_\lambda.$
The {\bf affine $q$-Schur algebra}, $\Sb_D$,
is the endomorphism ring of the right $\Hb_D$-module $\Tb_D.$ 
For any $\pen\in\Pen_D$ let ${}^\sharp\pen\in\NN^n$ be the $n$-uple
such that ${}^\sharp\pen_i=\sharp\pen^{-1}(i)$ for all $i\in[1,n]$.
Given $\lambda,\mu\in\Cen_D$ put
$${\ts\Sen_{\lambda\mu}=\{\sen\in\Sen_D\,|\,
\sum_{j\in\ZZ}\sen_{ij}={}^\sharp\lambda_i,\,
\sum_{i\in\ZZ}\sen_{ij}={}^\sharp\mu_j\}.}$$
An element $\sen\in\Sen_{\lambda\mu}$ is identified with the 
class in $S_\lambda\setminus S_D/S_\mu$ of the elements $w$ 
such that $(L_{(\lambda)w},L_\mu)\in Y_\sen$. Let 
$\Hb_{\lambda\mu}\subseteq\Hb_D$ be the $\Ab$-linear span
of the elements $T_\sen$ with $\sen\in\Sen_{\lambda\mu},$
where $T_\sen=\sum_{w\in\sen}T_w$. 
There is an isomorphism $\bigoplus_{\lambda\mu}\Hb_{\lambda\mu}\to\Sb_D$
which maps $T_\sen$, $\sen\in\Sen_{\lambda\mu}$, to the endomorphism such that
$T_\nu\mapsto\delta_{\mu\nu}T_\sen.$ Here $T_\sen$ is viewed as an element
in $T_\lambda\Hb_D$. 
It is well-known that $\CC_{G_D}[X_D\times X_D]$ is isomorphic to the 
specialization of $\Hb_D$ at $v=1/q$. Moreover we have the 
following result (see \cite {VV2, Proposition 7.4}).

\vskip3mm

\noindent{\bf Proposition.} {\it 
(a) The map $\Sb_{D|v=1/q}\to\CC_{G_D}[Y_D\times Y_D]$ such that
$T_\sen\mapsto\1b_\sen$ is an isomorphism of algebras.\hfill\break
(b) The map $\Tb_{D|v=1/q}\to\CC_{G_D}[Y_D\times X_D]$
such that $T_\lambda\mapsto\1b_\lambda$ 
extends uniquely to an isomorphism 
of $\bigl(\Sb_D\times\Hb_D\bigr)_{|v=1/q}$-modules.}
\qed

\vskip3mm

\noindent For any $\sen\in\Sen_D$ and $\pen\in\Pen_D,$ put 
$[\sen]=v^{y_\sen}T_\sen\in\Sb_D$ and $[\pen]=v^{x_\pen}T_\pen\in\Tb_D.$
Consider the $\CC$-linear involution on $\Hb_D$ such that 
$\bar T_w=T^{-1}_{w^{-1}},$ for all $w\in S_D$, and $\bar v=v^{-1}.$
It is well-known that $\bar{[\lambda]}=[\lambda]$ and that 
$\bar\Hb_{\lambda\mu}=\Hb_{\lambda\mu}$ for all $\lambda,\mu$. 
Let $\tau$ denote the antilinear involutions on $\Sb_D$, $\Tb_D$, such that
$$\tau([\pen])=\bar{[\pen]},\qquad\tau([\sen])=v^{-2x_\mu}\bar{[\sen]},\quad
\forall\sen\in\Sen_{\lambda\mu},\quad\forall\pen\in\Pen_D.$$ 
We have the following relations between the involutions on $\Sb_D$, $\Tb_D$,
and $\Hb_D$ (see \cite{VV2, Lemma 8.4}) : 
for all $y\in\Sb_D$, $x\in\Tb_D$, $z\in\Hb_D$, 
$$\tau(yx)=\tau(y)\tau(x)\and\tau(xz)=\tau(x)\bar z.\leqno(2.1)$$

\vskip3mm

\noindent{\bf 2.4.}
For any non-negative integer $k$ put $[k]=(v^k-v^{-k})/(v-v^{-1}),$ and 
$[k]!=[k][k-1]\cdots[1].$
The {\bf modified algebra} of $\widehat{\gen\len}_n$ at level 0
is the $\Kb$-algebra without unity $\Umod$
generated by the elements $\eb_i\ab_\mu$ and $\fb_i\ab_\mu,$ 
with $i\in[0,n-1]$ and $\mu\in\NN^n$, 
where the elements $\eb_i,\fb_i,\ab_\mu$
satisfy the following relations
$$\ab_\nu\ab_\mu=\ab_\mu\ab_\nu=\delta_{\nu\mu}\ab_\mu,\quad
\ab_{\mu+\omega_i}\eb_i=\eb_i\ab_{\mu+\omega_{i+1}},\quad
\ab_{\mu+\omega_{i+1}}\fb_i=\fb_i\ab_{\mu+\omega_i},$$
$$[\eb_i,\fb_j]\ab_\mu=\delta_{ij}[\mu_i-\mu_{i-1}]\ab_\mu,$$
$$\sum_{k=0}^{1-a_{ij}}(-1)^k\eb_i^{(k)}\eb_j\eb_i^{(1-a_{ij}-k)}=
\sum_{k=0}^{1-a_{ij}}(-1)^k\fb_i^{(k)}\fb_j\fb_i^{(1-a_{ij}-k)}=0\quad
\text{if}\quad i\neq j,$$
where $\eb_i^{(k)}=\eb_i^k/[k]!$, $\fb_i^{(k)}=\fb_i^k/[k]!,$
$(a_{ij})_{ij}$ is the Cartan matrix of type $A^{(1)}_{n-1}$, and
the $n$-uple $\omega_i$ is such that $\omega_{i,j}=\delta_{ij}$.
As usual, let also $\Ub^\pm$ be the $\Kb$-algebra generated by $\eb_i$ 
(resp. $\fb_i$), $i\in[0,n-1]$, modulo the Serre relations above.
Let $\tau$ be the involution of the $\CC$-algebras $\Ub^\pm$, $\Umod$, such that
$$\tau(v)=v^{-1},\qquad\tau(\eb_i)=\eb_i,\qquad\tau(\fb_i)=\fb_i,\qquad
\tau(\ab_\mu)=\ab_\mu.$$
For all $\lambda\in\Cen_D$ let ${}^\delta\lambda\in\Sen_D$ be the matrix 
such that ${}^\delta\lambda_{kl}=\delta_{kl}{}^\sharp\lambda_k$ for all $k,l$.
The following affine analogue of \cite{BLM} was observed independentely
by Lusztig and Ginzburg-Vasserot (see \cite{L1}, \cite{GV}, and also 
\cite{VV1,2}).

\vskip3mm

\noindent{\bf Proposition.} {\it There is  
a unique algebras homomorphism $\Phi_D\,:\,\Umod\to\Sb_D\otimes_\Ab\Kb$ such 
that $\Phi_D(\ab_{{}^\sharp\lambda}\eb_i)=[\sen],$ 
$\Phi_D(\fb_i\ab_{{}^\sharp\lambda})=[\ten],$
$\Phi_D(\ab_{\mu})=0$ if $\sum_i\mu_i\neq D$, 
and $\Phi_D(\ab_{{}^\sharp\lambda})=[{}^\delta\lambda],$ where 
$$\sen_{kl}=\ten_{lk}=
{}^\delta\lambda_{kl}-\delta_{ki}\delta_{li}+\delta_{ki}\delta_{l,i+1}\qquad
\forall k,l\in\ZZ.\leqno(2.2)$$
Moreover, we have $\Phi_D\circ\tau=\tau\circ\Phi_D$.
}\qed

\vskip3mm

\noindent{\bf Remark.} $(a)$ Observe that $y_\sen={}^\sharp\lambda_i-1$, 
$y_\ten={}^\sharp\lambda_{i+1}$, and $y_{{}^\delta\lambda}=0$.\hfill\break
$(b)$ Observe moreover that the map 
$\Cen_D\to\{\mu\in\NN^n\,|\,\sum_i\mu_i=D\}$, 
$\lambda\mapsto{}^\sharp\lambda$, is a bijection.

\vskip3mm

\noindent As a consequence, $\Tb_D\otimes_\Ab\Kb$ is a module over 
the modified algebra $\Umod$.

\vskip3mm

\noindent{\bf Proposition.} {\it The following formulas hold :
$$\matrix
{\ds \eb_i([\pen])=\sum_{k\in\pen^{-1}(i+1)}
v^{\sharp\{l\in\pen^{-1}(i+1)\,|\,l>k\}-\sharp\{l\in\pen^{-1}(i)\,|\,l>k\}}
[\pen_k^-],}\hfill\cr\cr
{\ds \fb_i([\pen])=\sum_{k\in\pen^{-1}(i)}
v^{\sharp\{l\in\pen^{-1}(i)\,|\,l<k\}-\sharp\{l\in\pen^{-1}(i+1)\,|\,l<k\}}
[\pen_k^+],}\hfill
\endmatrix\leqno(2.3)$$
where $\pen^\pm_k\in\Pen_D$ is the function such that
$$\pen_k^\pm(k)=\pen(k)\pm 1\quad\text{and}\quad
\pen_k^\pm(i)=\pen(i)\quad\text{if}\quad i\neq k\quad\mod D.$$
}

\vskip3mm

\noindent{\it Proof:} By \cite{VV1, Sections 5,6}, we have
$$\matrix
{\ds \eb_i(\1b_\pen)=v^{{}^\sharp\pen_i}\sum_{k\in\pen^{-1}(i+1)}
v^{-2\sharp\{l\in\pen^{-1}(i)\,|\,l>k\}}\1b_{\pen_k^-},}\hfill\cr\cr
{\ds \fb_i(\1b_\pen)=v^{{}^\sharp\pen_{i+1}}\sum_{k\in\pen^{-1}(i)}
v^{-2\sharp\{l\in\pen^{-1}(i+1)\,|\,l<k\}}\1b_{\pen_k^+}.}\hfill
\endmatrix$$
Thus it suffices to observe that
$$\matrix
x_\pen-x_{\pen_k^-}=\sharp\{l\st l>k,\,\pen(l)=\pen(k)\}-
\sharp\{l\st l<k,\,\pen(l)=\pen(k)-1\},\hfill\cr\cr
x_\pen-x_{\pen_k^+}=\sharp\{l\st l<k,\,\pen(l)=\pen(k)\}-
\sharp\{l\st l>k,\,\pen(l)=\pen(k)+1\}.\hfill\cr\cr
\endmatrix$$
\qed

\vskip1cm

\centerline{\smc 3. The crystal graph of $\Tb_D$.}

\vskip3mm

\noindent{\bf 3.1.}
Formulas (2.3) show that the $\Umod$-module $\Tb_D\otimes_\Ab\Kb$ is integrable.
Let $\tilde\eb_i$, $\tilde\fb_i$, be the Kashiwara operators
(see \cite{K1}). 
Let $\Lb_D$ be the $\Zb$-submodule of $\Tb_D$ linearly spanned by the
elements $[\pen]$. For each $\pen$ let $\b_\pen$ be the class of
$[\pen]$ in $\Lb_D/v\Lb_D$, and set $\B_D=\{\b_\pen\st\pen\in\Pen_D\}$.
Recall that the couple $(\Lb_D,\B_D)$ is a (lower) {\bf crystal base} if 
and only if the following properties hold :

\vskip2mm

\itemitem{-} for all $i$ we have 
$\tilde\eb_i(\Lb_D),\tilde\fb_i(\Lb_D)\subseteq\Lb_D,$

\vskip2mm

\itemitem{-} for all $i$ we have 
$\tilde\eb_i(\B_D),\tilde\fb_i(\B_D)\subseteq\B_D\cup\{0\},$

\vskip2mm

\itemitem{-} for all $i$  and all $\b,\b'\in\B_D$ we have 
$\tilde\eb_i(\b)=\b'$ if and only if $\tilde\fb_i(\b')=\b.$

\vskip3mm

\noindent{\bf Theorem.} {\it The couple $(\Lb_D,\B_D)$ is a crystal
base for the integrable $\Umod$-module $\Tb_D\otimes_\Ab\Kb$.}

\vskip3mm

\noindent{\it Proof:} The proof is very similar to the proof of
\cite{MM, Theorem 3.2.}. Fix $i\in[0,n-1]$. For any $\pen\in\Pen_D$ we 
consider a partition of the set $\pen^{-1}\bigl(\{i,i+1\}\bigr)$ 
into disjoints subsets $J,K_1,K_2,...,K_t,$ such that 

\vskip2mm

\itemitem{$(a)_{\ }$}we have $\sharp K_1=\sharp K_2=\cdots=\sharp K_t=2$, 

\vskip2mm

\itemitem{$(b)_\pen$} we have $\bigl(\pen(k),\pen(l)\bigr)\neq (i,i+1)$ if 
$k,l,$ are consecutive integers in $J$,

\vskip2mm

\itemitem{$(c)_\pen$} if $K_s=\{k,l\}$, $k<l$, $s\in[1,t]$,
then $\bigl(\pen(k),\pen(l)\bigr)=(i,i+1)$ and $[k,l]\cap J=\emptyset$.

\vskip2mm

\noindent Such a partition exists and is unique, up to the numbering of the 
pairs $K_1,K_2,...,K_t$. Given subsets $A\subseteq J$ and $B\subseteq[1,t]$,
let $\pen_{_{A,B}}\in\Pen_D$ be the function such that
$$\matrix
\pen^{-1}_{_{A,B}}(l)=\pen^{-1}(l)\hfill\quad&\text{if}\  l\neq i,i+1\,\mod\, n,
\hfill\cr\cr
\pen_{_{A,B}}(k)=i+1\hfill\quad&\forall k\in A,\hfill\cr\cr
\pen_{_{A,B}}(k)=i\hfill\quad&\forall k\in J-A,\hfill\cr\cr
\bigl(\pen_{_{A,B}}(k),\pen_{_{A,B}}(l)\bigr)=(i+1,i)\hfill\quad&
\text{if}\quad\{k,l\}=K_s,\quad k<l,\quad s\in B,\hfill\cr\cr
\bigl(\pen_{_{A,B}}(k),\pen_{_{A,B}}(l)\bigr)=(i,i+1)\hfill\quad&
\text{if}\quad\{k,l\}=K_s,\quad k<l,\quad s\notin B.\hfill
\endmatrix$$
Then, put
$$\langle\pen\rangle=\sum_{(A,B)}v^{n_A}(-v)^{\sharp B}[\pen_{_{A,B}}],$$
where the sum is taken over the set of couples $(A,B)$ such that 
$$A\subseteq J,\qquad B\subseteq[1,t],\qquad
\sharp A=\sharp\bigl(J\cap\pen^{-1}(i+1)\bigr),$$
and
$n_A=\sharp\{(k,l)\st k>l,\,k\in A,\,l\in J-A\}.$
Then, formulas (2.3) give
$$\eb_i(\langle\pen\rangle)=\sum_{A\subseteq J}\sum_{k\in A}
\sum_{B\subseteq[1,t]}v^{m_{A,k}+n_{A-k}}(-v)^{\sharp B}
[\pen_{_{A-k,B}}],$$
where
$$m_{A,k}=\sharp\{j\in J-A\st j<k\}-\sharp\{j\in J-A\st j>k\}.$$
Suppose now that $\pen^{-1}(i+1)=\emptyset$ and fix 
a partition $J,K_1,K_2,...,K_t$ of $\pen^{-1}(i)$ satisfying $(a)$.
For any integer $l\in[0,\sharp J]$, let $\pen_l\in\Pen_D$ be the unique
function such that 
$$\matrix
\pen^{-1}_l(j)=\pen^{-1}(j)\ \text{if}\  j\neq i,i+1\,\mod\, n,\hfill\cr\cr
J,K_1,K_2,...,K_t\ \text{satisfies}\ (b)_{\pen_l},\ (c)_{\pen_l},\hfill\cr\cr
\sharp\bigl(J\cap\pen_l^{-1}(i+1)\bigr)=l.\hfill
\endmatrix$$
Obviously, the map $(\pen, J, K_1,..., K_t, l)\mapsto \pen_l$
is a bijection onto $\Pen_D$. 
For any $A\subseteq J$ we have
$$\sum_{k\in J-A}v^{m_{A,k}}=[\sharp J-\sharp A].$$ 
Thus we get 
$$\eb_i(\langle\pen_l\rangle)=[\sharp J-l+1]\langle\pen_{l-1}\rangle
\and\fb_i(\langle\pen_l\rangle)=[l+1]\langle\pen_{l+1}\rangle,$$
where $\langle\pen_{-1}\rangle$ and $\langle\pen_{\sharp J+1}\rangle$ are
zero by definition. Now it suffices to observe that for any 
$\pen\in\Pen_D$ we have $\langle\pen\rangle\in[\pen]+v\Lb_D.$
\qed

\vskip3mm

\noindent{\bf Remark.} Observe that \cite{L1} implies the weaker statement
that $\pm\B_D$ is a signed crystal base.

\vskip3mm

\noindent{\bf 3.2.} 
Let $\Bmod\subset\Umod$ be the (global) canonical base.
Let $\Lmod\subset\Umod$ be the $\Zb$-submodule spanned by the elements
of $\Bmod$, and let $\dot\B$ be the projection of $\Bmod$ into $\Lmod/v\Lmod$.
Similarly, let $\Bb(\infty)\subset\Ub^-$ be the (global) canonical base,
let $\Lb(\infty)\subset\Ub^-$ be the $\Zb$-submodule spanned by the elements
of $\Bb(\infty)$, and let $\B(\infty)$ be the projection of $\Bb(\infty)$ 
into $\Lb(\infty)/v\Lb(\infty)$.
For any $\lambda\in\Cen_D$, $\pen\in\Pen_\lambda$, put
$$\bb_\pen=\sum_{\qen}\sum_{i\in\ZZ}v^{-i+x_\pen-x_\qen}
\dim\Hc^i_{X_{\qen,\lambda}}(IC_{X_{\pen,\lambda}})\,[\qen]\in\Lb_D,$$
where $X_{\pen,\lambda}$ is the fiber of the first projection 
$X_\pen\to Y_D$ at $L_\lambda$. 
By definition we have $\bb_\pen=[\pen]$ modulo $v\Lb_D.$
Thus the set $\Bb_D=\{\bb_\pen\st\pen\in\Pen_D\}$
is a $\Zb$-base of $\Lb_D$.
Recall that the affine $q$-Schur algebra $\Sb_D$ is 
identified with
the direct sum $\bigoplus_{\lambda,\mu\in\Cen_D}\Hb_{\lambda\mu}$. 
For any $\sen\in\Sen_{\lambda\mu}$ set 
$$\bb_\sen=\sum_{\ten}\sum_{i\in\ZZ}v^{-i+y_\sen-y_\ten}
\dim\Hc^i_{Y_{\ten,\lambda}}(IC_{Y_{\sen,\lambda}})\,[\ten]
\in\Hb_{\lambda\mu},$$
where $Y_{\sen,\lambda}$ is the fiber of the first projection 
$Y_\sen\to Y_D$ at $L_\lambda$. 
The set $\Bmod_D=\{\bb_\sen\st\sen\in\Sen_D\}$ is a $\Ab$-base of $\Sb_D$.
Let $\Lmod_D$ be the $\Zb$-linear span of $\Bmod_D$, and let $\dot\B_D$ be the
projection of $\Bmod_D$ in $\Lmod_D/v\Lmod_D$.

\vskip3mm

\noindent{\bf Proposition.} {\it Fix $\lambda\in\Cen_D$.
The following properties hold.\hfill\break
(a) We have $\Bb(\infty)\bigl([\lambda]\bigr)\subseteq\{0\}\cup\Bb_D$.
Moreover, if $\bb\bigl([\lambda]\bigr)=\bb'\bigl([\lambda]\bigr)$
with $\bb\neq\bb'\in\Bb(\infty)$, then $\bb\bigl([\lambda]\bigr)=0$.\hfill\break
(b) We have $\Lmod\bigl([\lambda]\bigr)\subseteq\Lb_D$.\hfill\break
(c) We have $\dot\B\bigl([\lambda]\bigr)\subseteq\{0\}\cup\B_D$.
Moreover, if $\b\bigl([\lambda]\bigr)=\b'\bigl([\lambda]\bigr)$
with $\b\neq\b'\in\dot\B$, then $\b\bigl([\lambda]\bigr)=0$.\hfill
}

\vskip3mm

\noindent{\it Proof:}
Claim $(a)$ follows from \cite{VV2, Remark 7.6}.
Let $\Qb$ be the generic Hall algebra of finite dimensional
nilpotent representations 
of the cyclic quiver of type $A^{(1)}_{n-1}$.
By definition, $\Qb$ is an $\Ab$-algebra.
It is well-known that $\Ub^-$ is a subalgebra of $\Qb\otimes_\Ab\Kb$. 
By \cite{VV2, Proposition 7.6}
there is an algebras homomorphism $\Theta\,:\,\Qb\to\Sb_D$
such that $\Theta$ is compatible with the involutions $\tau$ on $\Ub^-$ 
and $\Sb_D$, $\Theta\bigl(\Lb(\infty)\bigr)\subseteq\Lmod_D$,
and $\Theta\bigl(\B(\infty)\bigr)[{}^\delta\lambda]\subseteq\{0\}\cup\dot\B_D,$ 
for all $\lambda\in\Cen_D$.
Thus $\Theta\bigl(\Bb(\infty)\bigr)[{}^\delta\lambda]\subseteq\{0\}\cup\Bmod_D$.
By the definition of $\Theta$ (see [VV2]) we have
$$\Theta(\fb_i)[{}^\delta\lambda]=\Phi_D(\fb_i\ab_{{}^\sharp\lambda}),
\qquad\forall i.$$
Recall that, by the definition of the action of $\Sb_D$ on $\Tb_D$, we have
$$x(T_\lambda)=x=x[{}^\delta\lambda],
\qquad\forall x\in\Hb_{\mu\lambda}\subset\Sb_D.$$
Hence,
$$\Bb(\infty)\bigl([\lambda]\bigr)=
\Theta\bigl(\Bb(\infty)\bigr)\bigl([\lambda]\bigr)=
v^{x_\lambda}\Theta\bigl(\Bb(\infty)\bigr)[{}^\delta\lambda].$$
Moreover, for each $\mu,$ we have
$v^{x_\lambda}\Bmod_D\cap\Hb_{\mu\lambda}=\Bb_D\cap\Hb_{\mu\lambda}.$
The second part of Claim $(a)$ follows from the fact that
(see \cite{VV2, Proposition 7.6})
if $\Theta(\b)[{}^\delta\lambda]=\Theta(\b')[{}^\delta\lambda]$ and 
$\b\neq\b'\in\B(\infty)$, then $\Theta(\b)[{}^\delta\lambda]=0$.
Claims $(b)$ and $(c)$ follow from the following lemma.

\vskip3mm

\noindent{\bf Lemma.} 
{\it Let $M$ be an integrable $\Umod$-module and let $(\Lb_M,\B_M)$ be a 
crystal base of $M$. Let $\Psi\,:\,\Umod\ab_\lambda\to M$ be a $\Umod$-linear
homomorphism.\hfill\break
(a) If $\Psi\bigl(\Lb(\infty)\ab_\lambda\bigr)\subseteq\Lb_M$ 
and $\Psi\bigl(\B(\infty)\ab_\lambda\bigr)\subseteq\{0\}\cup\B_M$, then 
$\Psi\bigl(\Lmod\ab_\lambda\bigr)\subseteq\Lb_M$ and 
$\Psi\bigl(\dot\B\ab_\lambda\bigr)\subseteq\{0\}\cup\B_M$.\hfill\break
(b) Suppose moreover that $\Psi(\b\ab_\lambda)=0$ whenever 
$\Psi(\b\ab_\lambda)=\Psi(\b'\ab_\lambda)$ with $\b\neq\b'\in\B(\infty)$. 
Then the same property holds with $\b,\b'\in\dot\B$.
}

\vskip3mm

\noindent{\it Proof:}
Claim $(a)$ is already proved in \cite{K2, Proposition 9.1.3}.
More precisely, $\Psi$ splits through 
$\psi\,:\,\Umod\ab_\lambda\to V(\xi)\otimes\Lambda(\eta)$, $\ab_\lambda\mapsto
u_\xi\otimes u_\eta$, where $\xi,-\eta$, are high enough dominant weights with
$\xi-\eta=\lambda$, $V(\xi)$ and $\Lambda(\eta)$ are the simple modules
with highest weight $\xi$ and lowest weight $\eta$, and
$u_\xi\in V(\xi)$ (resp. $u_\eta\in\Lambda(\eta)$) is the highest weight vector 
(resp. the lowest weight vector).
Let $\B(\xi),\B(\eta),$ be the crystal bases of $V(\xi),\Lambda(\eta)$.
Claim $(a)$ is a consequence of the following property (see \cite{K2}) :
$$\matrix
\forall\b\in\B(\xi)\otimes\B(\eta),
\quad\exists\, i_1,i_2,...,i_N\in [0,n-1],\quad
\exists\, a_1,a_2,...,a_N\in \NN^\times\cr\cr
\text{such\ that}\quad
\tilde\fb^{(a_N)}_{i_N}\cdots\tilde\fb^{(a_2)}_{i_2}\tilde\fb_{i_1}^{(a_1)}
(\b)\in\B(\xi)\otimes u_\eta.\endmatrix
\leqno(3.0)$$
From (3.0) and the fact that 
$\tilde\fb_i\bigl(\B(\xi)\otimes u_\eta\bigr)\subseteq
\{0\}\cup\B(\xi)\otimes u_\eta$ we get :
if $\b\neq\b'\in\B(\xi)\otimes\B(\eta)$ then
there are $i_1,i_2,...,i_N\in [0,n-1]$ and
$a_1,a_2,...,a_N\in \NN^\times$ such that
$$ \tilde\fb^{(a_N)}_{i_N}\cdots\tilde\fb^{(a_2)}_{i_2}\tilde\fb_{i_1}^{(a_1)}
(\b)=\b_1\otimes u_\eta\and
\tilde\fb^{(a_N)}_{i_N}\cdots\tilde\fb^{(a_2)}_{i_2}\tilde\fb_{i_1}^{(a_1)}
(\b')=\b'_1\otimes u_\eta,$$
where $\b_1\neq\b_1'\in\B(\xi)\cup\{0\}$.
Thus, Claim $(b)$ results from the (known) fact that if
$\b\neq\b'\in\dot\B$, then $\psi(\b\ab_\lambda)\neq\psi(\b'\ab_\lambda)$,
unless they are both equal to zero.
\qed

\vskip3mm

\noindent{\bf 3.3.}
By definition of the involution $\tau$ on $\Tb_D$ we have 
$\tau(\bb_\pen)=\bb_\pen$, for all $\pen\in\Pen_D.$

\vskip3mm

\noindent{\bf Proposition.} {\it For all $\bb\in\Bmod$ and all 
$\lambda\in\Cen_D$ we have $\bb\bigl([\lambda]\bigr)\in\{0\}\cup\Bb_D.$
Moreover, if $\bb\bigl([\lambda]\bigr)=\bb'\bigl([\lambda]\bigr)$
with $\bb\neq\bb'\in\Bmod$, then $\bb\bigl([\lambda]\bigr)=0$.
}

\vskip3mm

\noindent{\it Proof:}
By Proposition 3.2 we have $\dot\B\bigl([\lambda]\bigr)\subseteq\{0\}\cup\B_D$.
Now, by (2.1) the elements in $\Bmod\bigl([\lambda]\bigr)$
are stable by the involution $\tau$.
Thus, $\Bmod\bigl([\lambda]\bigr)\subseteq\{0\}\cup\Bb_D.$ 
The second claim follows from the second claim of Proposition 3.2.$(c)$.
\qed

\vskip1cm

\centerline{\smc 4. The compatibility of $\Phi_D$ with the bases.} 

\vskip3mm

\noindent{\bf 4.1.} 
The following theorem is the main result of the paper.

\vskip3mm

\noindent{\bf Theorem.} {\it For all $\bb\in\Bmod$ 
we have $\Phi_D(\bb)\in\{0\}\cup\Bmod_D.$ Moreover, the kernel of
$\Phi_D$ is linearly spanned by the elements $\bb\in\Bmod$ such that
$\Phi_D(\bb)=0$.}

\vskip3mm

\noindent{\it Proof:}
For any $\lambda,\mu\in\Sen_D$, let $\Phi_{\lambda\mu}$ be the component 
of $\Phi_D$ in $\Hb_{\lambda\mu}$.
By construction, for any $x\in\Umod$ the element 
$\Phi_{\lambda\mu}(x)$ is the component of $v^{-x_\mu}x\bigl([\mu]\bigr)$
in $\Hb_{\lambda\mu}\subseteq T_\lambda\Hb_D$. 
It is known that $\Bmod=\bigcup_{\lambda\mu}\ab_\lambda\Bmod\ab_\mu$.
Moreover, we have $\Phi_D(x)=\Phi_{\lambda\mu}(x)$ for all 
$x\in\ab_\lambda\Umod\ab_\mu$.
Thus the theorem follows from Proposition 3.3, since 
$v^{-x_\mu}\Bb_D\cap\Hb_{\lambda\mu}=\Bmod_D\cap\Hb_{\lambda\mu}$ by definition.
\qed

\vskip3mm

\noindent{\bf 4.2.}
Fix positive integers $D_1,D_2$ such that $D=D_1+D_2$.
Set $g=-Id_{\K^{D_1}}+Id_{\K^{D_2}}\in G_D.$
The fixpoints set $Y_D^g$ is isomorphic to $Y_{D_1}\times Y_{D_2}$.
More precisely, if $\lambda\in\Cen_D$ let $Y_\lambda$ be the $G_D$-orbit
of $L_\lambda$. Then,  
$$Y_\lambda^g=\bigcup_{\lambda_1,\lambda_2}Y_{\lambda_1}\times Y_{\lambda_2},$$
where the sum is taken over the couples $(\lambda_1,\lambda_2)$ such that
${{}^\sharp\lambda_1+{}^\sharp\lambda_2={}^\sharp\lambda}.$
Put $L_D=G_{D_1}\times G_{D_2}$ and let $U_D$ be the unipotent radical of the 
stabilizer of the flag $\{0\}\subseteq \K^{D_1}\subseteq \K^{D}.$ 
For any element $L\in Y^g_D$ the $U_D$-orbit of $L$ is a 
(possibly infinite dimensional) $\k$-linear space.
It may be viewed as the set of flags $\tilde L\in Y_D$ such that 
$\exp(\var\,g)(\tilde L)\to L$ when $\var\to\infty$. 
We will use the following notation :
$$\tilde L\to L\quad\iff\quad\tilde L\in U_D(L).$$
For any $\sen\in\Sen_D$ and any couple $(L,L')\in Y_D^g\times Y_D^g$, 
the set
$$\{(L,\tilde L')\in Y_\sen\,|\,\tilde L'\to L'\}$$
is finite.
Consider the map 
$\Omega_{D_1D_2}\,:\,\CC_{G_D}[Y_D\times Y_D]\to\CC_{L_D}[Y^g_D\times Y^g_D]$
such that
$$\Omega_{D_1D_2}(\alpha)(L,L')=
\sum_{\tilde L'\to L'}\alpha(L,\tilde L').\leqno(4.1)$$
Then 
$$\Omega_{D_1D_2}(\alpha\star\beta)(L,L")
=\sum_{\tilde L"\to L"}\alpha\star\beta\,(L,\tilde L")
=\sum_{\tilde L'\in Y_D}\alpha(L,\tilde L')\sum_{\tilde L"\to L"}
\beta(\tilde L',\tilde L").$$
Given $\tilde L'\in Y_D$ there is an element $u\in U_D$ such
that $u\bigl(\tilde L'\bigr)\in Y^g_D$. 
Since $\beta$ is $G_D$-invariant, the second sum does not depend  
on $\tilde L'$ but only on the limit of $\exp(\var\,g)(\tilde L')$
when $\var\to\infty$. 
Thus $\Omega_{D_1D_2}(\alpha\star\beta)=
\Omega_{D_1D_2}(\alpha)\star\Omega_{D_1D_2}(\beta).$

\vskip3mm

\noindent Let $\Delta\,:\,\Umod\to\Umod\otimes\Umod$ be the $\Kb$-linear map
such that
$$\matrix
{\ds\Delta(\ab_\lambda)=\sum_{\lambda_1+\lambda_2=\lambda}
\ab_{\lambda_1}\otimes\ab_{\lambda_2},}\hfill\cr\cr
{\ds\Delta(\ab_\lambda\eb_i)=\sum_{\lambda_1+\lambda_2=\lambda}
\ab_{\lambda_1}\otimes\ab_{\lambda_2}\bigl(
v^{\lambda_{1,i}}\otimes\eb_i+\eb_i\otimes v^{-\lambda_{2,i}}
\bigr),}\hfill\cr\cr
{\ds\Delta(\fb_i\ab_\lambda)=\sum_{\lambda_1+\lambda_2=\lambda}
\ab_{\lambda_1}\otimes\ab_{\lambda_2}\bigl(
v^{-\lambda_{1,i+1}}\otimes\fb_i+\fb_i\otimes v^{\lambda_{2,i+1}}
\bigr).}\hfill
\endmatrix$$

\vskip3mm

\noindent{\bf Proposition.} {\it The map $\Omega_{D_1D_2}$ is the specialization
at $v=1/q$ of an algebras homomorphism $\Sb_{D}\to\Sb_{D_1}\otimes\Sb_{D_2}$,
still denoted by $\Omega_{D_1D_2}$. Moreover we have
$$\Omega_{D_1D_2}\circ\Phi_D=\bigl(\Phi_{D_1}\otimes\Phi_{D_2}\bigr)\circ\Delta.
\leqno(4.2)$$
}

\vskip3mm

\noindent{\it Proof:} Fix $\lambda\in\Cen_D$.
Let us check the equality (4.2) on the generators $\ab_{{}^\sharp\lambda}$,
$\ab_{{}^\sharp\lambda}\eb_i$ and $\fb_i\ab_{{}^\sharp\lambda}$. 
If $\alpha\in\CC_{G_D}[Y_\lambda\times Y_\lambda]$ is the characteristic 
function of the diagonal, then only $\tilde L'=L'$ contributes to (4.1), and
we get the equality (4.2) for $\ab_{{}^\sharp\lambda}.$ 
We consider now the case of the element $\ab_{{}^\sharp\lambda}\eb_i$.
Thus, suppose that ${}^\sharp\lambda-\omega_i+\omega_{i+1}\in\NN^n$.
Fix $\lambda_1,\lambda_2,$ such that 
${}^\sharp\lambda_1+{}^\sharp\lambda_2={}^\sharp\lambda$.
Consider the matrix $\sen\in\Sen_D$ such that
$\sen_{kl}={}^\delta\lambda_{1,kl}+{}^\delta\lambda_{2,kl}-
\delta_{ki}\delta_{li}+\delta_{ki}\delta_{l,i+1}$ (see (2.2)). 
Fix moreover $L=(L^1,L^2)\in Y_{\lambda_1}\times Y_{\lambda_2}$
and $L'=({L'}^1,{L'}^2)\in Y_{D_1}\times Y_{D_2}$.
Let $E(L,L')$ be the set of the flags $\tilde L'\in Y_D$
such that $(L,\tilde L')\in Y_\sen$ and $\tilde L'\to L'$. Then,
$$E(L,L')=\{(V_j)_j\in Y_D\,|\,
L_{j-1}\subseteq V_j\subseteq L_j,\,\dim(L_j/V_j)=\delta_{ij},\,$$
$${L'}^1_j=V_j\cap\K^{D_1},\,{L'}^2_j=\rho(V_j)\},$$
where $\rho\,:\,\K^D\to\K^{D_2}$ is the projection along $\K^{D_1}$.
In particular, if the set $E(L,L')$ is non empty then 
$L^a_{i-1}\subseteq {L'}^a_i\subseteq L^a_i$ for $a=1,2,$
and one of the following two cases holds. 

\vskip2mm

\itemitem{-} Either ${L'}^1_i=L^1_i$ and 
$L^2_{i-1}\subseteq{L'}^2_i\subsetneq L^2_i$.
Then, $E(L,L')=\{L^1_i\oplus{L'}^2_i\}.$

\vskip2mm

\itemitem{-} Either ${L'}^2_i=L^2_i$ and 
$L^1_{i-1}\subseteq{L'}^1_i\subsetneq L^1_i$.
Then,
$$E(L,L')=\{V\,|\,{L'}^1_i\oplus L^2_{i-1}\subset V,\,
L^1_i\oplus V=L_i\}\simeq\k^{{}^\sharp\lambda_{2,i}}.$$

\noindent Hence,
$$\Omega_{D_1D_2}\bigl([\sen]\bigr)=
v^{y_\sen-y_{\sen-{}^\delta\lambda_1}}[{}^\delta\lambda_1]\otimes
[\sen-{}^\delta\lambda_1]
+v^{y_\sen-y_{\sen-{}^\delta\lambda_2}-2{}^\sharp\lambda_{2,i}}
[\sen-{}^\delta\lambda_2]\otimes[{}^\delta\lambda_2].$$
The result follows from the identities (see Remark 2.4)
$$y_\sen-y_{\sen-{}^\delta\lambda_1}={}^\sharp\lambda_{1,i}
\and
y_\sen-y_{\sen-{}^\delta\lambda_2}-2{}^\sharp\lambda_{2,i}=-{}^\sharp\lambda_{2,i}.$$
The case of $\fb_i\ab_{{}^\sharp\lambda}$ is identical. 
\qed

\vskip3mm

\noindent Recall that $\Sb_n=\bigoplus_{\lambda\mu}\Hb_{\lambda\mu}$ where
$\lambda,\mu\in\Pen_n$. Moreover $\Hb_{\lambda\mu}=\Hb_n$ if
$\lambda=\mu=(1,2,...,n).$
There is an algebras homorphism $\eps\,:\,\Sb_n\to\Ab$ which is zero on 
$\Hb_{\lambda\mu}$ if $\lambda$ or $\mu$ is not $(1,2,...,n)$, and such that
$\eps(T_i)=-1$ for all $i$ (view each component $\Hb_{\lambda\mu}$ as a 
subspace in $\Hb_D$ and  make it acts on the sign representation of $\Hb_D$). 
Let $\phi\,:\,\Umod\to\Umod$ be the endomorphism such that 
$$\phi(\ab_\lambda)=\ab_{\lambda-(1,1,...,1)},\qquad
\phi(\eb_i)=v\,\eb_i,\and\phi(\fb_i)=v^{-1}\fb_i,$$
where $\ab_\lambda$ is zero if $\lambda\notin\NN^n$.
The map $\Phi'_{D+n,D}=(\eps\otimes 1)\Omega_{nD}$
is an algebras homomorphism $\Sb_{D+n}\to\Sb_D$ such that
$\Phi'_{D+n,D}\circ\Phi_{D+n}=\Phi_D\circ\phi$, by $(4.2)$.
Let $\psi$ be the automorphism of $\Sb_D$ which acts on the component
$\Hb_{\lambda\mu}$ be the scalar $v^{\sum_i(\lambda_i-\mu_i)}$.
Put $\Phi_{D+n,D}=\psi\circ\Phi'_{D+n,D}.$
The map $\Phi_{D+n,D}$ is precisely the algebras homomorphism
$\Sb_{D+n}\to\Sb_D$ introduced in \cite{L1, Section 9.1}.
In particular, we have
$$\Phi_{D+n,D}\circ\Phi_{D+n}=\Phi_D.\leqno(4.4)$$
Let $\leq$ be the standard order on $\Sen_D$, i.e. $\ten\leq\sen$ if and
only if $Y_\ten\subseteq\bar Y_\sen.$ In particular, if
$\ten\leq\sen$ then $\ten_{ii}\geq\sen_{ii}$ for all $i$.

\vskip3mm

\noindent{\bf Lemma.} {\it Fix $\sen\in\Sen_{D+n}$ and let $\ten$
be the $\ZZ\times\ZZ$-matrix such that $\ten_{ij}=\sen_{ij}-\delta_{ij}$ 
for all $i,j\in\ZZ$. If $\ten_{ij}\geq 0$ for all $i,j$, then
$\Phi_{D+n,D}\bigl([\sen]\bigr)=c\,[\ten]$ modulo some lower term,
where $c$ is a non-zero constant.
}

\vskip3mm

\noindent{\it Proof:}
Fix a couple $(L,L')$ in $Y^g_D\times Y^g_D$.
If the couple $(L,L')$ belongs to the support of 
$\Omega_{D_1D_2}\bigl([\sen]\bigr)$,
then there is an element $(L,\tilde L')$ in $Y_\sen$ such that
$\exp(\eps g)(L,\tilde L')\to (L,L')$ when $\eps\to\infty$.
Thus, $(L,L')\in\bar Y_\sen$.
Put $L=(L_1,L_2)$, $L'=(L'_1,L'_2)$, where $L_i,L'_i\in Y_{D_i}$. 
Fix $\ten_i\in\Sen_{D_i}$ such that $(L_i,L'_i)\in Y_{\ten_i}$. 
Let $\ten=\ten_1+\ten_2\in\Sen_D$ be the matrix such that
$\ten_{ij}=\ten_{1,ij}+\ten_{2,ij}$ for all $i,j$.
Then, $(L,L')\in Y_{\ten}$.
Thus if $Y_{\ten_1}\times Y_{\ten_2}$ is in the support of
$\Omega_{D_1D_2}\bigl([\sen]\bigr)$, for any matrices
$\ten_i\in\Sen_{D_i}$ and $\sen\in\Sen_D$, then $\ten_1+\ten_2\leq\sen$.
In particular, let us consider the case $D\to D+n$, $D_1\to n$, $D_2\to D$,
and fix $\sen\in\Sen_D$ such that $\sen_{ii}\geq 1$ for all $i$.
If the orbit $Y_{\ten}$, $\ten\in\Sen_D$, is in the support of
$\Phi_{D+n,D}\bigl([\sen]\bigr)$ then 
there is a matrix $\ten_1\in\Sen_n$ such that $\ten_1+\ten\leq\sen.$
Moreover since $\eps$ is zero on $\Hb_{\lambda\mu}$ if $\lambda$ or $\mu$ 
is not $(1,2,...,n)$, necessarily $\ten_{1,ii}\leq 1$ for all $i$. 
Hence, $\ten_{ii}\geq\sen_{ii}-1$ for all $i$.
Among all such orbits, the bigger one is labelled by
the matrix $\ten_0$ with $\ten_{0,ij}=\sen_{ij}-\delta_{ij}$.
Put $\lambda=(1,2,...,n)$. For any matrix $\ten_1\in\Sen_n$
such that $Y_{\ten_1}\times Y_{\ten_0}$ is in the support of
$\Omega_{nD}\bigl([\sen]\bigr)$, we have $\ten_1+\ten_0\leq\sen$
and, thus, $\ten_{1,ii}\geq 1$ for all $i$. If moreover 
$\ten_1\in\Sen_{\lambda\lambda}$, then $\ten_{1,ij}=\delta_{ij}$ for all $i,j$.
Hence, to prove that $Y_{\ten_0}$ is indeed in the support of
$\Phi_{D+n,D}\bigl([\sen]\bigr)$ it suffices to prove that 
$Y_{Id}\times Y_{\ten_0}$ is in the support of
$\Omega_{nD}\bigl([\sen]\bigr)$. This is obvious since $Id+\ten_0=\sen$. 
\qed

\vskip3mm

\noindent{\bf Remark.} The map $\Omega_{D_1D_2}$ was given in the non-affine 
case in \cite{G}. After this preprint was written, Lusztig gave us a copy
of \cite{L2} where \cite{L2, Conjecture 9.2} is proved for the case 
$\gen\len_2$. The preprint \cite{L2} contains also a construction of the map 
$\Phi_{D+n,D}$.

\vskip3mm

\noindent{\bf 4.3.} 
Following \cite{L1}, let $\Sen^{ap}_D\subset\Sen_D$ be the set of matrices 
$\sen$ such that for any $j\in\ZZ-\{0\}$ there exist $i\in\ZZ$ with 
$\sen_{i,i+j}=0$. The matrices in $\Sen^{ap}_D$ are said to be 
{\bf aperiodic}. By \cite{L1, Theorem 8.2}, the subfamilly
$\Bmod_D^{ap}=\{\bb_\sen\,|\,\sen\in\Sen^{ap}_D\}\subset\Bmod_D$ 
is a $\Kb$-basis of $\Im\Phi_D$. Given an element $\bb\in\Bmod^{ap}_{D+n}$,
by Theorem 4.1 there is an element $\bb'\in\Bmod$ such that
$\Phi_{D+n}(\bb')=\bb$. Then, by (4.4), we have
$\Phi_{D+n,D}(\bb)=\Phi_D(\bb')$. Hence, Theorem 4.1 implies that 
$\Phi_{D+n,D}(\Bmod^{ap}_{D+n})\subseteq\{0\}\cup\Bmod^{ap}_D$.
Now, we have $\Im\Phi_D=\Phi_{D+n,D}\bigl(\Im\Phi_{D+n}\bigr)$ by (4.4).
Thus, $\Phi_{D+n,D}(\Bmod^{ap}_{D+n})=\{0\}\cup\Bmod^{ap}_D$.
Moreover, formula $(4.4)$ and the second claim in Theorem 4.1 guarantee that
if two distinct elements $\bb,\bb'\in\Bmod^{ap}_{D+n}$ have the same image
by $\Phi_{D+n,D}$, then this image is $0$.
The following more precise result was conjectured in
\cite{L1, Conjecture 9.2}.

\vskip3mm

\noindent{\bf Theorem.} {\it Fix $\sen\in\Sen^{ap}_{D+n}$ and let $\ten$
be the $\ZZ\times\ZZ$-matrix such that $\ten_{ij}=\sen_{ij}-\delta_{ij}$ 
for all $i,j\in\ZZ$. Then\hfill\break
(a) we have $\Phi_{D+n,D}(\bb_\sen)=0$ if 
$\ten_{ij}<0$ for some $i,j$,\hfill\break
(b) we have $\Phi_{D+n,D}(\bb_\sen)=\bb_\ten$
if $\ten_{ij}\geq 0$ for all $i,j$.
}

\vskip3mm

\noindent{\it Proof:}
We know that $\Phi_{D+n,D}(\bb_\sen)=0$ or $\bb_\uen$ for some $\uen$.
Suppose that $\ten_{ij}=\sen_{ij}-\delta_{ij}$.
If $\ten_{ij}\geq 0$ for all $i,j$, then by Lemma 4.2 we have
$\Phi_{D+n,D}(\bb_\sen)=\bb_\ten$. Claim $(b)$ is proved.
Let $\Bmod^{ap\,*}_{D+n}\subset\Bmod^{ap}_{D+n}$ be the subfamilly
labelled by the matrices $\sen$ with $\sen_{ii}\geq 1$ for all $i$.
Claim $(b)$ implies that $\Phi_{D+n,D}(\Bmod^{ap\,*}_{D+n})=\Bmod^{ap}_D$.
Thus, Claim $(a)$ follows from 
$\Phi_{D+n,D}(\Bmod^{ap}_{D+n})=\{0\}\cup\Bmod^{ap}_D$.
\qed

\vskip1cm

\noindent{\it Acknowledgements.}
{\eightpoint{Part of this work was done while the second author was visiting
the Institute for Advanced Study at Princeton. The second author
is grateful to G. Lusztig for his kind invitation.}}

\vskip3cm
{\eightpoint{
$$\matrix\format\l&\l&\l&\l\\
\phantom{.} & {\text{Olivier Schiffmann}}\phantom{xxxxxxxxxxxxx} &
{\text{Eric Vasserot}}\\
\phantom{.}&{\text{D\'epartement de Math\'ematiques}}\phantom{xxxxxxxxxxxxx} &
{\text{D\'epartement de Math\'ematiques}}\\
\phantom{.}&
{\text{Ecole Normale Sup\'erieure}}
\phantom{xxxxxxxxxxxxx} &
{\text{Universit\'e de Cergy-Pontoise}}
\\
\phantom{.}&
{\text{45 rue d'Ulm}}
\phantom{xxxxxxxxxxxxx} & 
{\text{2 Av. A. Chauvin}}
\\
\phantom{.}&
{\text{75005 Paris}}
\phantom{xxxxxxxxxxxxx} & 
{\text{95302 Cergy-Pontoise Cedex}}
\\
\phantom{.}&{\text{France}}\phantom{xxxxxxxxxxxxx} & 
{\roman{France}}\\
&{\text{email: schiffma\@clipper.ens.fr}}\phantom{xxxxxxxxxxxxx} &
{\text{email: vasserot\@math.pst.u-cergy.fr}}
\endmatrix$$
}}

\enddocument